\documentclass[letterpaper, 10 pt, journal, twoside]{ieeecolor}
\usepackage{lcsys}
\usepackage[left=0.7in,top=0.7in,right=0.7in,bottom=0.7in]{geometry}
\usepackage{mwe}
\usepackage{afterpage}
\usepackage{graphics} 
\usepackage{epsfig,url} 
\usepackage{times} 
\usepackage{bigints}

\usepackage{amsmath} 
\usepackage{amssymb,xcolor}  
\usepackage{mathtools,gensymb,amsfonts,amsthm,mathrsfs,comment,bm}
\usepackage{booktabs}
\usepackage[capitalise]{cleveref}
\usepackage{color}

\theoremstyle{plain}
\newtheorem{theorem}{Theorem}
\newtheorem{lemma}[theorem]{Lemma}

\theoremstyle{definition}

\theoremstyle{remark}
\newtheorem{remark}{Remark}

\newcommand{\N}{\mathbb{N}}

\DeclarePairedDelimiter\abs{\lvert}{\rvert}
\DeclarePairedDelimiter\norm{\lVert}{\rVert}
\makeatletter
\let\oldabs\abs
\def\abs{\@ifstar{\oldabs}{\oldabs*}}
\let\oldnorm\norm
\def\norm{\@ifstar{\oldnorm}{\oldnorm*}}
\makeatother

\begin{document}

\title{The Exponential Stabilization of a Heat and Piezoelectric Beam Interaction with Static or Hybrid Feedback Controllers}
\author{Ahmet \"Ozkan \"Ozer \IEEEmembership{Member, IEEE}, Ibrahim Khalilullah, Uthman Rasaq
\thanks{A.\"O. \"Ozer is with Department of Mathematics, Western Kentucky University,
	Bowling Green, KY 42101, USA.  (e-mail: ozkan.ozer@wku.edu).}
\thanks{I. Khalilullah and U. Rasaq are graduate students at the Department of Mathematics, Western Kentucky University (WKU),
	Bowling Green, KY 42101, USA. (e-mails: uthman.rasaq114@topper.wku.edu, skmdibrahim.khalilullah504@topper.wku.edu).}
	\thanks{As part of A.\"{O}. \"{O}zer's sabbatical leave in 2023,  Enrique Zuazua's generous support is appreciated at the hosting institute Friedrich-Alexander-Universität Erlangen-Nürnberg in Germany. Moreover, the support from the NSF of USA under Cooperative Agreement No. 1849213 and the RCAP grant of WKU are acknowledged.}
}

\maketitle
\thispagestyle{empty}

\begin{abstract}
This study investigates a strongly-coupled system of partial differential equations (PDE) governing heat transfer in a copper rod, longitudinal vibrations, and total charge accumulation at electrodes within a magnetizable piezoelectric beam. Conducted within the transmission line framework, the analysis reveals profound interactions between traveling electromagnetic and mechanical waves in magnetizable piezoelectric beams, despite disparities in their velocities. Findings suggest that in the open-loop scenario, the interaction of heat and beam dynamics lacks exponential stability solely considering thermal effects. To confront this challenge, two types of boundary-type state feedback controllers are proposed: (i) employing static feedback controllers entirely and (ii) adopting a hybrid approach wherein the electrical controller dynamically enhances system dynamics. In both cases, solutions of the PDE systems demonstrate exponential stability through meticulously formulated Lyapunov functions with diverse multipliers. The proposed proof technique establishes a robust foundation for demonstrating the exponential stability of Finite-Difference-based model reductions as the discretization parameter approaches zero.
\end{abstract}


\section{Introduction}

Piezoelectric materials, including potassium sodium niobate, barium titanate, or bismuth sodium titanate, possess the unique ability to generate electric displacement directly proportional to an applied mechanical stress \cite{Yang}. These materials are integral in the development of actuators, where electrical input—typically voltage—is essential \cite{O-IEEE2,AMOP}. The drive frequency, a crucial aspect of the electrical input, governs the speed at which a piezoelectric beam vibrates or transitions between states. Moreover, piezoelectric materials play roles as sensors and energy harvesters in various applications \cite{Baur,Kiran}.

While electrostatic approximations based on Maxwell’s equations suffice for describing low-frequency vibrations in many piezoelectric applications, electromagnetic effects are either fully or partially disregarded. However, for devices like piezoelectric acoustic wave devices, magnetic effects may play a more substantial role, tightly coupled to mechanical effects. Existing models in the literature, relying solely on electrostatic approximations, fall short in capturing the vibrational dynamics of these devices. Hence, there is a critical need for more accurate models that consider both electromagnetic and mechanical interactions, as demonstrated in previous works such as \cite{D,T} and related references. Employing the complete set of Maxwell equations in the modeling process enables the appropriate coupling of electromagnetic waves to mechanical vibrations, a concept termed piezo-electro-magnetism \cite{Yang}.

The incorporation of a fully dynamic theory, especially in modeling (acoustic) magnetizable piezoelectric beams, induces a significant shift in controllability dynamics compared to models obtained through electrostatic approximation. While the electrostatic model is precisely controllable by a tip-velocity controller \cite{M-O,AMOP}, the fully dynamic model faces challenges, particularly the lack of controllability of high-frequency vibrational modes \cite{M-O,Voss,AMOP} with the same controller design. Consequently, adapting the controller design becomes imperative to achieve exact controllability in the context of the fully dynamic model.

The underlying partial differential equation (PDE) model for this physics problem shares similarities with the fluid-structure interactions found in the linearized one-dimensional Navier-Stokes equations (heat equation) and the solid structure (piezoelectric beam), see e.g. \cite{Zua1}. These interactions hold paramount importance in various scenarios, encompassing aerodynamic considerations in aircraft airflow \cite{Castille} and the application of piezoelectric energy harvesters in underwater environments \cite{Kiran}.

\subsection{System of PDEs  under Consideration}
A copper beam of length \(l_1\), clamped on the left end and free on the other end, is attached to the piezoelectric beam of length \(l_2\), free at both ends. Piezoelectric beams are elastic smart beams with electrodes at their top and bottom surfaces, connected to an external electric circuit.
Two controllers are considered for the piezoelectric beam at its right end, one for controlling the mechanical strains and another one for controlling the total current accumulated at the electrodes. The heat on the copper beam transmits through the right end, and the heat and its flux at the joint interact directly with the mechanical strains on the beam. With the small displacements assumption due to the linear Euler-Bernoulli beam theory, transverse oscillations of the piezoelectric beam are relatively negligible in comparison to the longitudinal vibrations in the form of expansion and compression of the centerline of the beam.
Noting that \(\alpha, \alpha_1, \gamma, \beta, \mu, \rho\) are the piezoelectric material parameters, define
\begin{eqnarray}
\label{matrices1}
\begin{array}{ll}
 M:=\begin{bmatrix}
    \rho &0  \\
  0& \mu        \\
\end{bmatrix} , A:=\begin{bmatrix}
  \alpha       & -\gamma \beta   \\
   -\gamma\beta       & \beta  \\
\end{bmatrix}, \alpha=\alpha_1+\gamma^2\beta,
\end{array}
\end{eqnarray}
where the stiffness coefficient \(\alpha\), due to the involvement of piezoelectricity, is different from that of fully elastic materials \(\alpha_1\). Denote \(z(x,t), v(x,t),\) and \(p(x,t)\) by the heat distribution on the copper rod, longitudinal oscillations of the centerline of the beam, and the total charge accumulated at the electrodes of the beam, respectively. The equations of motion are a system of strongly-coupled partial differential equations,
\begin{eqnarray}
		\left\{
		\label{main}  \begin{array}{ll}
	z_t(x,t)-\kappa z_{xx}(x,t)=0, & x \in (-l_1,0),\\
			\begingroup 
M
\endgroup  \begin{bmatrix}
				{v_{tt}}  \\
				{p}_{tt}       \\
			\end{bmatrix}(x,t)-
A \begin{bmatrix}
				v_{xx}  \\
				p_{xx}       \\
			\end{bmatrix}(x,t)=0, & x\in (0,l_2),\\
		z(-l_1,t)= p(0,t)=0, \\
		z(0,t) = v_t(0,t),& \\
		\kappa z_x(0,t)=\alpha v_x(0,t)-\gamma\beta p_x(0,t),&\\
A
 \begin{bmatrix}
				v_{x}  \\
				p_{x}   \\
			\end{bmatrix}(l_2,t) = \begin{bmatrix}
				g_1 \\
				g_2   \\
			\end{bmatrix}(t),&  t\in \mathbb{R}^+,\\
			z(x,0)=z^0(x), & x\in [-l_1,0],\\
\left(v, p,\dot v, \dot p \right)(x,0)\\
\qquad =(v^0,p^0, v^1,p^1)(x),&x\in [0,l_2]
		\end{array}\right.
	\end{eqnarray}
where  $\kappa$ is the thermal diffusivity constant, $\rho$, $\alpha$ (and $\alpha_1$), $\beta$, $\gamma$, and $\mu$ are all piezoelectric material-specific positive constants, namely, the mass density per unit volume, the elastic and piezoelectric stiffness, the beam coefficient of impermeability, the piezoelectric constant, and the magnetic permeability, $g_1(t), g_2(t)$ are the strain and voltage controllers.

\begin{figure}[htb!]
		\centering
		{{\includegraphics[width=7.5  cm]{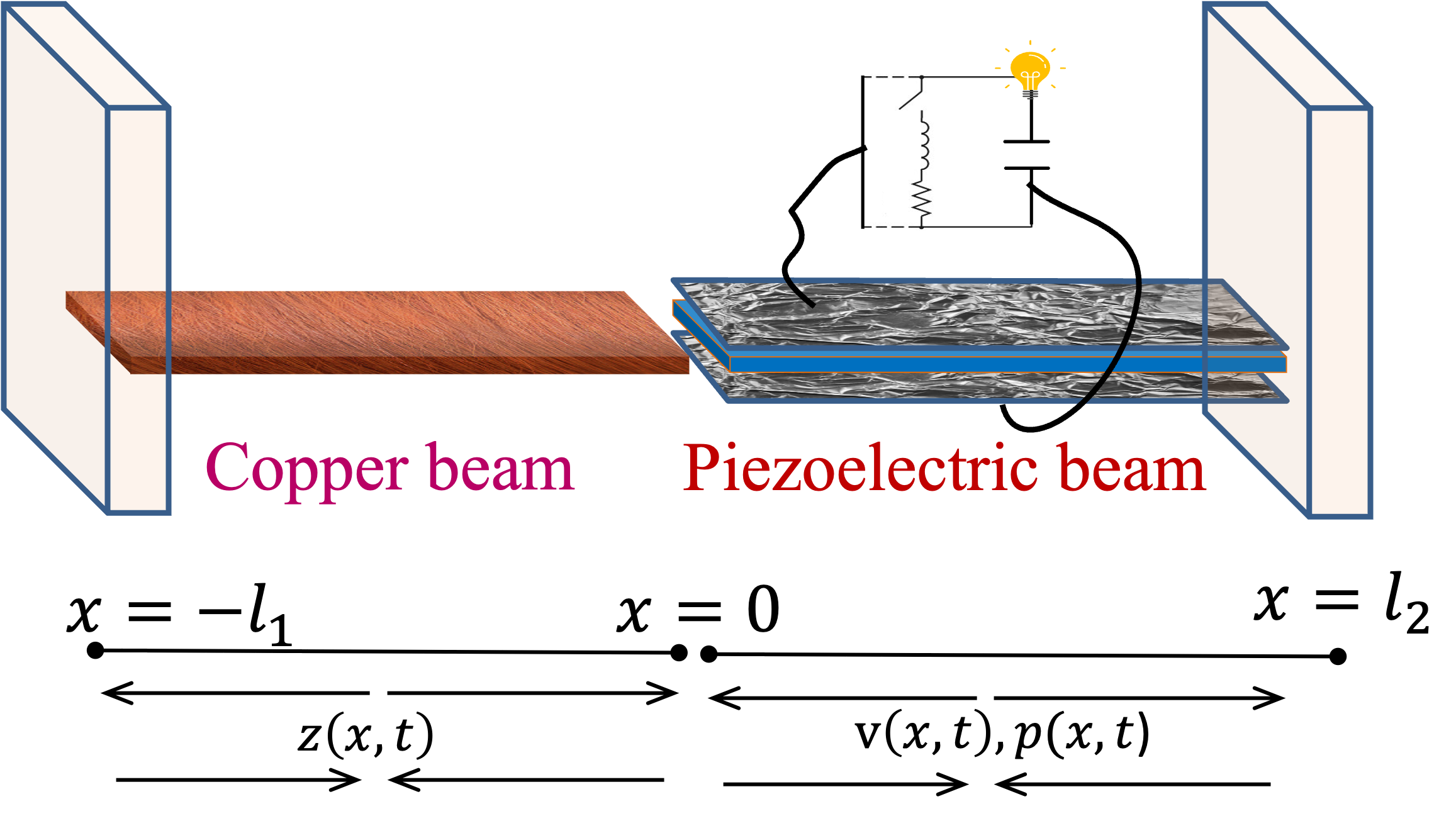} }}%
		\vspace{-0.1in}
\caption{\footnotesize
Magnetic effects contribute to beam vibrations and electric field across electrodes. Simultaneously, heat transmits from the left copper beam to the right piezoelectric beam. Continuity conditions for $z(\cdot,t)$ \& $v_t(\cdot,t)$ and $z_x(\cdot,t)$ \& $\alpha v_x (\cdot,t)-\gamma\beta p_x (\cdot,t)$ are maintained at joint $x=0$. System dynamics are exponentially stabilized by two boundary feedback controllers.}
		\label{single}%
	\end{figure}

The energy of the system \eqref{main} is defined by
\begin{eqnarray}\label{eq4}
\begin{array}{ll}
E(t)=\frac{1}{2}\int_0^{l_1}\underbrace{\left|z(-x,t)\right|^2}_{\text {Thermal~ energy}} dx
       \\
       +\frac{1}{2}\int_0^{l_2}\left[\underbrace{\rho \left|v_t\right|^2}_{\text {Mechanical ~kinetic~ energy}}+\underbrace{ \mu  \left|p_t\right|^2}_{\text {Magnetic~ kinetic~ energy}}\right.
       \\\left.+\underbrace{\alpha_1 |v_x|^2}_{\text {Potential~ energy}}+\underbrace{\beta |\gamma v_x-p_x|^2}_{\text {Electromechanical ~ energy}}\right] (x,t)dx, ~t\ge0.
       \end{array}
\end{eqnarray}
Note that it is known that with only one feedback controller, the piezoelectric beam dynamics cannot be exponentially stabilized for almost all material parameters \cite{M-O,O-MCSS,Voss}. It is also known that the exponential stability of the heat equation above is not enough to exponentially stabilize a serially-connected wave equation \cite{Lassy,Zua1} without a boundary feedback controller \cite{Zheng}. In this model, two feedback controllers are extremely needed \cite{O-R,Ramos}.
\subsection{First-order Formulation}
\label{Sec2}
Define  spatial derivatives by $D_x^k=\frac{\partial^k}{\partial x^k},~ k=1,2,\ldots.$ Let $\vec\psi=[u_1,u_2,w_1,w_2]^{\rm T},$
Let $u_1:=v_x(x,t)$ $u_2:=p_x(x,t),$ $w_1:=v_t(x,t),$ $w_2:=p_t(x,t),$  and
$\tilde M_{4\times 4} :=\begin{bmatrix}
		I& 0 \\
			0& M      \\
			\end{bmatrix},$ $\tilde A_{4\times 4}:=\begin{bmatrix}
		0 & I \bigotimes D_x\\
			A\bigotimes D_x& 0      \\
			\end{bmatrix}$
where $\bigotimes$ is the Kronecker product. Then, \eqref{main} can be rewritten as
\begin{eqnarray}
\label{main2}
		\left\{
		\label{main2}  \begin{array}{ll}
	z_t(x,t)-\kappa D_x^2 z (x,t)=0, & x \in (-l_1,0),\\
\tilde M\dot {\vec \psi}(x,t)-
\tilde A\vec \psi(x,t)=0, & x\in (0,l_2),\\
		z(-l_1,t)=w_2(0,t)=0, \\
		z(0,t) = w_1(0,t),  \\
		\kappa z_x(0,t)=\alpha u_1 (0,t)-\gamma\beta u_2(0,t),&\\
A
 \begin{bmatrix}
				u_1 \\
				u_2   \\
			\end{bmatrix}(l_2,t) = \begin{bmatrix}
				g_1 \\
				g_2   \\
			\end{bmatrix}(t),&  t\in \mathbb{R}^+,\\
			z(x,0)=z^0(x), & x\in [-l_1,0]\\
\left(u_1,u_2,w_1, w_2 \right)(x,0),\\
\qquad =(v^0_x, p^0_x, v^1,p^1)(x),&x\in [0,l_2].
		\end{array}\right.
	\end{eqnarray}
Define the natural energy space as $\mathrm H=L^2(0,l_1)\times \left[L^2(0,l_2)\right]^4,$
equipped with the inner product
\begin{eqnarray*}
\label{inner}
\begin{array}{ll}
\left<\vec \Phi,{\vec{\tilde \Phi}}\right>_{\mathrm H}=\frac{1}{2}\int_0^{l_1} \Phi_1(-x) {\bar{\tilde \Phi}}_1(-x)  dx
       \\
       +\frac{1}{2}\int_0^{l_2}\left[\rho \Phi_4(x) {\bar{\tilde \Phi}}_4(x)+\mu \Phi_5(x) {\bar{\tilde \Phi}}_5(x)+\alpha_1 \Phi_{2,x}(x){\bar{\tilde \Phi}}_2(x) \right. \\
       \left.+\beta \left(\gamma\Phi_{2,x}(x)-\Phi_{3,x}(x)\right) \left(\gamma {\bar{\tilde\Phi}}_{2,x}(x)- {\bar{\tilde\Phi}}_{3,x}(x)\right) \right] dx,
\end{array}
\end{eqnarray*}
for  $t\in \mathbb{R}^+$ and $\vec \Phi,{\vec{\tilde \Phi}}\in\mathrm{H}.$  With the choice of the state vector $\vec X(t)=[z(-x,t),u_1(x,t),u_2(x,t),w_1(x,t),w_2(x,t)]^{\rm T},$ \eqref{main2} can be formulated as
\begin{eqnarray*}
\label{inner}
\left\{
\begin{array}{ll}
{\dot {\vec X}}=\mathcal A \vec X:=\begin{bmatrix}
		I\bigotimes D_x^2& 0 \\
			0& (\tilde M^{-1} \tilde A)\bigotimes D_x       \\
			\end{bmatrix}\vec X,\\
 \vec X(0)=[z^0(-x),u^0_x(x),v^0_x(x),u^1(x),v^1(x)]^{\rm T}\in\mathrm{H}.
\end{array}\right.
\end{eqnarray*}
Here, the domain of $\mathcal A$ is defined by
	\begin{eqnarray*}
\label{dom2}
\left\{
\begin{array}{ll}{\rm Dom}(\mathcal A)=\left\{z, u_1,  u_2, v_1, v_2 \in \mathrm{H}: z(-x)\in H^2(0,l_1), \right.\\
u_1,  u_2, v_1, v_2\in H^1(0,l_2), z(-l_1,t)=0, \\
 \alpha u_1 (l_2,t) -\gamma\beta  u_2(l_2,t)=-\xi_1 v_1(l_2,t), \\
   \beta  u_2(l_2,t) -\gamma\beta  u_1(l_2,t)=-\xi_2 v_2(l_2,t), w_1(0,t)=z(0,t), \\
 \left.  w_2(0)=0,\quad   \alpha u_1(0,t) -\gamma\beta  u_2(0,t)=\kappa z_x(0,t) \right\}.
\end{array}\right.
\end{eqnarray*}
The following result is immediate. The proof is skipped here.
\begin{theorem} The operator $(\mathcal A, D(\mathcal A))$ generates a $C_0-$ semigroup on $\mathrm H.$
\end{theorem}
\subsection{Static and Hybrid Feedback Control Design}
\label{design}
	Two feedback control designs are studied in this paper.
\noindent{\bf Static Feedback Controllers:} This type of problem is more classical but well-appreciated due to its simplicity.
For $\xi_1,\xi_2>0$, chose the controllers in \eqref{main2} as the following
\begin{equation}\label{staticcon}\begin{bmatrix}
				g_1(t) \\
				g_2(t)   \\
			\end{bmatrix} =\begin{bmatrix}
				-\xi_1 v_t(l_2,t) \\
				-\xi_2 p_t (l_2,t)  \\
			\end{bmatrix}.
			\end{equation}
\noindent {\bf Hybrid (Static and Dynamic) Feedback Controllers:} Piezoelectric beams actuated by a dynamic voltage/charge source require a more elaborated controllers design. Therefore, letting $A\in \mathbb{R}^{n\times n} ,$  $\vec b, \vec c\in \mathbb{R}^n$ and $d,\xi_1>0,$  choose the controllers in \eqref{main2} as the following
\begin{eqnarray}\label{main-aug1}
\begin{array}{ll}
\begin{bmatrix}
				g_1(t) \\
				g_2 (t) \\
			\end{bmatrix} =\begin{bmatrix}
				-\xi_1 v_t(l_2,t) \\
				- \vec c^{\rm{T}} \vec q(t) -d p_t(l_2,t)
			\end{bmatrix},
			\end{array}
\end{eqnarray}
 where $\vec q(t)$ satisfies its own dynamic differential equation,
\begin{eqnarray}\label{main-aug2}
\begin{array}{ll}\vec q_t(t)= {\bf A} \vec q(t) + \vec b p_t (l_2,t),\quad \vec q (0)=\vec \zeta.
\end{array}
\end{eqnarray}
	For the hybrid case, the following are assumed  \cite{Morgul,Morgul2}:
\begin{itemize}
\item[1.] ${\bf A}$ has eigenvalues with negative real parts.
\item[2.] $({\bf A},\vec b)$ and $({\bf A},\vec c)$ are exactly controllable and exactly observable, respectively.
\item[3.] There exists a parameter $\Gamma\ge 0$ such that $d\ge \Gamma$ and $d+{\rm {Re}} ~\vec c^{\rm{T}} (isI-{\bf A})^{-1}\vec b>\Gamma.$
\end{itemize}
By the Meyer-Kalman-Yakubovich lemma, for any given $n\times n$ positive-definite matrix $Q$, there exist an $n\times n$ positive-definite matrix $P$,  an $n\times 1$ vector   $\vec q_1,$ and a constant $\Delta >0$ such that
\begin{eqnarray}\label{MKY}
\begin{array}{ll}
A^{\rm{T}}P+A P^{\rm{T}}=-\vec q_{1} \vec{q}_1^{\rm{T}}-\Delta Q,\\
\quad P \vec b-\vec c=\sqrt{2(d-\Gamma)} \vec q_1.
\end{array}
\end{eqnarray}
	\begin{remark} Note that for the voltage/charge-controlled electrostatic piezoelectric beam model, a scalar-type dynamic controller is proposed for the first time \cite{O-IEEE2,AMOP}. Indeed,
	 for $\xi_1,\xi_2>0$, the feedback controllers are chosen as
	$	\begin{bmatrix}
				g_1(t) \\
				g_2 (t) \\
			\end{bmatrix} =\begin{bmatrix}
				-\xi_1 v_t(l_2,t) \\
				-\P(t)   \\
			\end{bmatrix}
	$  where the charge $\P(t)$ satisfies
			\begin{eqnarray}\label{dumbb}
			 \P_t(t)= - \P(t)+ \xi_2\ p_t (l_2,t), \quad \P(0)=\eta.
			\end{eqnarray}  The proposed dynamic control law for the electrical part  extends the work in \cite{AMOP}.
	\end{remark}

	\subsection{Our Contribution}
	To the best of our knowledge, this is the first work devoted to the interaction of the heat equation with the fully-dynamic (magnetizable) piezoelectric beam in a transmission line setting. By designing static and hybrid feedback controllers, overall electromagnetic vibrations and heat distribution can be controlled to the equilibrium state exponentially fast. The  exponential stability result in each case is based on a carefully chosen Lyapunov function and multipliers, eliminating the spectral analysis completely, i.e. see \cite{Lassy,Jin,Zua1}. Hence, the recently proposed order-reduced model reduction by Finite Differences, as in [11], is applicable, which is one of the main objectives of the paper.

Both closed-loop systems make a solid foundation for the interaction of the linearized one-dimensional fluid equations (replacing the heat equation) and piezoelectric acoustic devices under the ocean, used for energy harvesting purposes \cite{Kiran}.
\section{Exponential Stability with Static Feedback Controllers}
\label{sec2}
\begin{lemma} \label{lem1} The system \eqref{main2}  with \eqref{staticcon} is dissipative, i.e.
\begin{eqnarray}\label{FD-exp5}
\begin{array}{ll}
\frac{dE(t)}{dt} &=-\xi_1 |w_1(l_2,t)|^2-\xi_2 |w_2(l_2,t)|^2\\
&-\int_0^{l_1}  \kappa |z_x(-x,t)|^2 dx \le 0.
\end{array}
\end{eqnarray}
\end{lemma}
\begin{proof}
By taking the time derivative of  $E(t)$ along the solutions of \ \eqref{main2}  with \eqref{staticcon}, and by the integration by parts on the respective domains,
\begin{eqnarray*}
\begin{array}{ll}
\frac{dE(t)}{dt}   =-\kappa z(-x,t)z_x(-x,t)|_0^{l_1}- \kappa\int_0^{l_1}{|{z}}_{x}(-x,t)|^2 dx\\
    \quad  + \left[w_1(x,t) (\alpha u_{1}(x,t) -\gamma\beta  u_{2}(x,t))+ w_{2}(x,t) (\beta u_{2}(x,t)\right.\\
    \left. -\gamma \beta  u_{1}(x,t))\right]|_0^{l_1}-\int_0^{l_2}\left[w_{1,x}(x,t) (\alpha u_{1}(x,t) -\gamma\beta  u_{2}(x,t))\right. \\
     \left. + w_{2,x}(x,t) (\beta u_{2}(x,t) -\gamma \beta  u_{1}(x,t))+ \alpha_1 (u_1 w_{1,x}) (x,t))\right.\\
     \left.+\beta (\gamma u_1(x,t)-u_2(x,t)) (\gamma w_{1,x}(x,t)- {w}_{2,x}(x,t))\right]dx.\\
\end{array}
\end{eqnarray*}
Recalling \eqref{matrices1} for $\alpha=\alpha_1+\gamma^2\beta$,  \eqref{FD-exp5} follows.
\end{proof}

For $\delta, a_1,b_1, c_1>0$, define the  Lyapunov functional
\begin{eqnarray}\label{Lyp1}
\begin{array}{ll}
L(t):=E(t)+ \delta \left(F_1(t)+F_2(t) +F_3(t)\right),\\
F_1(t):=a_1\int_0^{l_2} x \left\{\rho  u_1 w_1 + \mu  u_2 w_2\right\}(x,t)~dx, \\
F_2(t):= b_1 \int_0^{l_2} (l_2-x)\left[\alpha |u_1(x,t)|^2+ \beta |u_2(x,t)|^2 \right.\\
\left. -\gamma\beta (u_1u_2)(x,t)+\rho |w_1(x,t)|^2+\mu |w_2(x,t)|^2\right] dx,\\
 F_3(t):=c_1\int_0^{l_1} |z(-x,t)|^2 dx.
\end{array}
\end{eqnarray}
The following technical results are immediate.
\begin{lemma} \label{lem2} For the choices of $c_1:=b_1 l_2$ and
\begin{equation}\label{constants}
\begin{array}{ll}
a_1:=2b_1\left[\frac{4 l_2\kappa }{l_1^2} + \max\left(\sqrt{\frac{\alpha_1}{\rho}}+\frac{\gamma^2\sqrt{\beta\mu}}{\rho},2\sqrt{\frac{\beta}{\mu}}\right)\right],
\end{array}
\end{equation}
define the constant $M:=b_1 l_2 \tilde M$ where
\begin{eqnarray}\label{constant2}
\begin{array}{ll}
\tilde M:=\max\left\{ \left[\frac{8l_2\kappa }{l_1^2} +2 \max\left(\sqrt{\frac{\alpha_1}{\rho}}+\frac{\gamma^2\sqrt{\beta\mu}}{\rho},2\sqrt{\frac{\beta}{\mu}}\right)\right] \right.\\
\left.\max \left[\sqrt{\frac{\rho}{\alpha_1}}+ \sqrt{\frac{\mu\gamma^2}{\alpha_1}},\sqrt{\frac{\mu}{\beta}}+\sqrt{\frac{\mu\gamma^2}{\alpha_1}}, 2, \frac{\alpha}{\alpha_1}+\frac{\gamma^2\beta}{2\alpha_1}\right]\right\}.
\end{array}
\end{eqnarray}
Then, for $0<\delta < \frac{1}{M},$  $L(t)$ is equivalent to $E(t),$ i.e.,
\begin{equation}\label{eqrfd5}
(1-M \delta) E(t)\le L(t)\le (1+M\delta)E(t).
\end{equation}
\end{lemma}

\begin{proof}  By the Cauchy-Schwarz's, H\"older's, and Minkowski's inequalities,  $F_1(t)$ is estimated as the following
\begin{eqnarray*}
\begin{array}{ll}
F_1(t)\leq a_1l_2\int_0^{l_2}   \left\{ \sqrt\frac{\rho}{\alpha _1} \left(\frac{\alpha_1 }{2}|u_1|^2+\frac{\rho}{2} |w_1|^2\right)(x,t) \right.\\
\quad + \sqrt{\frac{\mu}{\beta}}  \left( \frac{\mu}{2} |w_2|^2   + \frac{\beta}{2}|u_2 -\gamma u_1|^2\right)(x,t) \\
 \quad \left. + \sqrt\frac{\mu\gamma^2}{\alpha_1}\left(\frac{\mu}{2} |w_2|^2 dx+ \frac{\alpha_1}{2}|u_1|^2\right)\right\}(x,t) dx\\
\quad \leq a_1 l_2 \max \left[\sqrt \frac{\rho}{\alpha_1}+\sqrt\frac{\mu\gamma^2}{\alpha_1},\sqrt{\frac{\mu}{\beta}}+\sqrt\frac{\mu\gamma^2}{\alpha_1}\right]E(t).
\end{array}
\end{eqnarray*}
Next, we estimate $F_2(t).$ By  \eqref{matrices1} and the Cauchy-Schwarz's inequality,
\begin{eqnarray*}
\begin{array}{ll}
&F_2(t)\leq b_1 l_2\int_0^{l_2}\left[\alpha_1 |u_1|^2+\beta|\gamma u_1-u_2|^2 +\rho |w_1|^2\right.\\
&~~\left.  +\frac{\gamma^2 \beta |u_1|^2}{2} +\frac{\beta|u_2-\gamma u_1+\gamma u_1|^2}{2}+\mu |w_2|^2\right](x,t) dx\\
& \leq 2\frac{b_1 l_2}{2}\int_0^{l_2}\left[\alpha_1 |u_1|^2+\beta|\gamma u_1-u_2|^2 +\rho |w_1|^2\right.\\
&\left.  +\frac{\gamma^2 \beta |u_1|^2}{2} +\beta|u_2-\gamma u_1|^2+\gamma ^2\beta|u_1|^2+\mu |w_2|^2\right](x,t) dx\\
&= 2\frac{b_1l_2}{2}\int_0^{l_2}\left\{\left(1+\frac{3\gamma^2 \beta}{2\alpha_1}\right)\alpha_1|u_1|^2+2\beta|\gamma u_1-u_2|^2\right.\\
&\quad\left.+\rho |w_1|^2+\mu |w_2|^2\right\}(x,t)dx\\
&\leq 2b_1l_2 \max \left(\frac{\alpha}{\alpha_1}+\frac{\gamma^2\beta}{2\alpha_1},2\right)E(t).
\end{array}
\end{eqnarray*}
Finally, $F_3(t)\le 2 c_1   E(t)$ simply follows.
By the choice constants  $a_1,c_1,$ and $\tilde M$ as in \eqref{constants}-\eqref{constant2}, \eqref{eqrfd5} follows.
\end{proof}

\begin{lemma} \label{lem3} Assume \eqref{constants}. Then,  $F(t)$ in \eqref{Lyp1} satisfies the following estimate
\begin{equation}\label{eqrfd2}
\begin{array}{ll}
 \frac{dF(t)}{dt}\le \frac{-8b_1l_2\kappa }{l_1^2} E(t)+\frac{a_1 l_2}{2} \left(\rho +\frac{2}{\alpha_1}\xi_1^2\right)|w_1(l_2,t)|^2\\
 ~+\frac{a_1 l_2}{2} \left(\mu +\frac{\alpha+\gamma^2\beta}{\alpha_1\beta}\xi_2^2\right)|w_2(l_2,t)|^2.
  \end{array}
\end{equation}
\end{lemma}
\begin{proof}
First, dot product the matrix equation of \eqref{main2}    by the vector of multipliers
\begin{equation}
\left(\begin{bmatrix} \label{giant}
		\rho x w_1&  (\alpha u_1 -\gamma \beta u_2)(l_2-x) \\
	\mu x w_2 &  \beta(u_2-\gamma u_1)(l_2-x)  \\
	 x u_1&   w_1 (l_2-x) \\
	x u_2  & w_2 (l_2-x)\\
			\end{bmatrix}(x,t) \right)\begin{bmatrix} a_1\\ b_1\end{bmatrix},
			\end{equation}and integrate the sum in $x$  from $0$ to $l_2$. After a couple of integration by parts, the boundary conditions \eqref{staticcon}, and Cauchy-Schwarz's inequality we have two identities below
\begin{eqnarray}
&\label{mult1}
\begin{array}{ll}
 a_1\frac{d}{dt}\int^{l_2}_0  x\left[\rho  w_1 u_1+\mu  u_2 w_2\right](x,t)dx\\
\quad  =-\frac{a_1}{2}\int^{l_2}_0 \left[\rho |w_1|^2 +\mu  |u_2|^2 \right.  \\
 \qquad\qquad\left.+\alpha_1 |u_1^2 + \beta |\gamma u_1-u_2|^2 \right](x,t)dx\\
\qquad +\frac{a_1 l_2}{2} \left\{ \rho |w_1|^2 +\mu |w_2|^2 +\frac{\left|\gamma \xi_2 w_2+ \xi_1 w_1\right|^2}{\alpha_1}\right. \\
\left. \qquad\qquad\qquad +\frac{\left|\xi_2 w_2\right|}{\beta}\right\}(l_2,t),
\end{array}\\
&\label{mult2_a}
\begin{array}{ll}
&\frac{b_1}{2}\frac{d}{dt}\int^{l_2}_0 (l_2 -x) \left[\alpha |u_1|^2 +\beta |u_2|^2\right.\\
&\quad \left.- \gamma \beta u_1 u_2 + \rho |w_1|^2 +\mu |w_2|^2\right] (x,t)dx\\
&\le -2 b_1 l_2 \kappa z_x(0,t)z(0,t)+2\int^{l_2}_0 \frac{b_1}{2} \left[\sqrt{\frac{\alpha_1}{\rho}} \alpha_1 |u_1|^2\right.\\
&\quad+\sqrt{\frac{\beta}{\mu}} \beta |\gamma u_1-u_2|^2 + 2\sqrt{\frac{\beta}{\mu}} \mu |w_2|^2\\
&\left. \quad+ \left(\sqrt{\frac{\alpha_1}{\rho}} +2\sqrt{\frac{\beta}{\mu}}\frac{\mu \gamma^2}{\rho}  \right) \rho |w_1|^2\right](x,t) dx.
\end{array}
\end{eqnarray}
Next, multiply both sides of the first  equation of \eqref{main2}  by $c_1 z(-x,t)$ integrate it in $x$  from $0$ to $l_2$ to get
\begin{eqnarray}\label{mult3}
\begin{array}{ll}
   &\frac{c_1}{2} \frac{d}{dt}\int^{l_1}_0 |z(-x,t)|^2 dx= -2\kappa c_1 z(0,t)z_x(0,t)\\
   &\qquad\qquad  -2\kappa c_1 \int_0^{l_1} |z_x(-x,t)|^2 dx.
\end{array}
\end{eqnarray}
By adding \eqref{mult2_a} to \eqref{mult3} together with the Wirtinger's inequality and  \eqref{constants},
\begin{eqnarray}
\label{mult23_a}
\begin{array}{ll}
&\frac{b_1}{2}\frac{d}{dt}\int^{l_2}_0 (l_2 -x)\left[\alpha |u_1|^2 +\beta |u_2|^2 - \gamma \beta u_1 u_2 \right.\\
&\left.+ \rho |w_1|^2+\mu |w_2|^2\right] (x,t) dx\\
&+ \frac{c_1}{2} \frac{d}{dt}\int^{l_1}_0 |z(-x,t)|^2 dx=2\int^{l_2}_0 \frac{b_1}{2} \left[\sqrt{\frac{\alpha_1}{\rho}} \alpha_1 |u_1|^2 \right.\\
&\left.+\sqrt{\frac{\beta}{\mu}} \beta |\gamma u_1-u_2|^2+ 2\sqrt{\frac{\beta}{\mu}} \mu |w_2|^2\right. \\
&\left. + \left(\sqrt{\frac{\alpha_1}{\rho}} +2\sqrt{\frac{\beta}{\mu}}\frac{\mu \gamma^2}{\rho}  \right) \rho |w_1|^2\right](x,t) dx\\
& -\frac{4b_1l_2 \kappa}{l_1^2} \int_0^{l_1} |z_x(-x,t)|^2 dx. \\
\end{array}
\end{eqnarray}
Finally, by adding \eqref{mult1} and \eqref{mult23_a}, and using  the Cauchy-Schwarz's inequality one more time,
\begin{eqnarray*}
\begin{array}{ll}
\frac{dF}{dt}&\le  \frac{-8b_1 l_2 \kappa}{l^2_1}\frac{1}{2}\int^{l_1}_0 |z(-x,t)|^2 dx\\
&~ +\frac{a_1 l_2}{2} \left(\rho +\frac{2}{\alpha_1}\xi_1^2\right)|w_1(l_2,t)|^2\\
&~ +\frac{a_1 l_2}{2} \left(\mu +\frac{\alpha+\gamma^2\beta}{\alpha_1\beta}\xi_2^2\right)|w_2(l_2,t)|^2\\
& ~-\frac{1}{2} \int_0^{l_2} \left[a_1\left\{\rho |w_1|^2+\mu |w_2|^2 \right. \right.\\
&\quad \left.+\alpha_1 |u_1|^2+\beta |\gamma u_1-u_2|^2 \right\}\\
&~+b_1 \left\{ 2\sqrt{\frac{\beta}{\mu}} \mu |w_2|^2+ \sqrt{\frac{\alpha_1}{\rho}} \alpha_1 |u_1|^2\right.\\
& \quad +\left(\sqrt{\frac{\alpha_1}{\rho}} + \frac{2\gamma^2}{\rho} \sqrt{\frac{\beta}{\mu}}\right) \rho|u_2|^2\\
&\left.\left.\quad +2\sqrt{\frac{\beta}{\mu}}|\gamma u_1-u_2|^2 \right\} \right](x,t)dx.
\end{array}
\end{eqnarray*}
Finally, \eqref{eqrfd2} follows from \eqref{constants}.
\end{proof}

  	\begin{theorem} \label{mainthm1}Suppose that there exists a positive constant
	\begin{equation*}\delta <{\rm {min}} \left(\frac{1}{M},\frac{2\xi_1}{a_1l_2 \left(\rho+\xi_1^2\frac{2}{\alpha_1}\right)},\frac{2\xi_2}{a_1l_2 \left(\mu+\xi_2^2\frac{\alpha+\gamma^2\beta}{\alpha_1\beta}\right)}\right).\end{equation*}
Then, the energy $E(t)$ along the solutions of  \eqref{main2}-\eqref{staticcon} decays exponentially with the decay rate $\sigma:=\delta (1-M\delta) \frac{8 b_1 l_2\kappa}{l_1^2},$
\begin{eqnarray*}
E(t)\le\frac{1+M\delta}{1-M\delta} e^{-\sigma t} E(0), \quad \forall t>0.
\end{eqnarray*}
\end{theorem}
\begin{proof} By a direct calculation with  Lemmas \ref{lem1}-\ref{lem3},
\begin{eqnarray*}
\begin{array}{ll}
\frac{dL}{dt}&\le - \left(\xi_1-\frac{\delta a_1 l_2}{2}\left(\rho +\frac{2}{\alpha_1}\xi_1^2\right)\right)|w_1(l_2,t)|^2\\
&\quad-\left(\xi_2 - \frac{\delta a_1 l_2}{2}\left(\mu +\frac{\alpha+\gamma^2\beta}{\alpha_1\beta}\xi_2^2\right)\right)|w_2(l_2,t)|^2\\
&\quad -\int_0^{l_1}  \kappa |z_x(-x,t)|^2 dx- \delta \frac{8 b_1 l_2\kappa }{l_1^2}E(t) \\
&\le - \delta \frac{8 b_1 l_2\kappa }{l_1^2}E(t).
\end{array}
\end{eqnarray*}
By Lemma \ref{lem2}, $\frac{dL}{dt}\le -\delta (1-M\delta) \frac{8 b_1 l_2 \kappa }{l_1^2} L(t).$
Next, by the direct application of  the Gr\"onwall's inequality and Lemma \ref{lem2} again,
\begin{eqnarray*}
\begin{array}{ll}
(1-M\delta) E(t)\le L(t)\le L(0)e^{-\sigma t} \le (1+M\delta) E(0) e^{-\sigma t}.
\end{array}
\end{eqnarray*}
Hence, the conclusion is immediate.
\end{proof}
Note that the decay rate $\sigma(\delta)$ is maximal as the system exponentially stabilizes fastest.  It is straightforward that $\sigma(\delta(\xi_1,\xi_2))$ makes its maximal value
$\sigma_{max}(\delta)= \frac{2b_1 l_2 \kappa  }{l_1^2 M} ~{\rm achieved ~at}~ \delta=\frac{1}{2M}.$
Since $\delta$ is a function of $\xi_1$ and $\xi_2$, to achieve $\delta=\frac{1}{2M},$ an optimization argument may be needed to provide safe intervals for the feedback amplifiers $\xi_1$ and $\xi_2,$ as in \cite{O-R}.

\section{Exponential Stability with Hybrid Feedback Controllers}	
\label{sec3}
The energy \eqref{eq4} of the new system is redefined as
\begin{eqnarray}\label{en-hyb}
\begin{array}{ll}
E_{\rm{hybrid}}(t)=E(t)+\frac{1}{2} \vec q^{\rm T} P \vec q, \quad t\in \mathbb{R}^+.
       \end{array}
\end{eqnarray}

    \begin{lemma}\label{lem6} The energy \eqref{en-hyb} for the closed-loop system \eqref{main2}  with \eqref{main-aug1}-\eqref{main-aug2}  is dissipative, i.e.
\begin{eqnarray}\label{dissp2}
\begin{array}{ll}
\frac{dE_{\rm{hybrid}}(t)}{dt}&=-\frac{1}{2}  \left[ \vec q^{\rm T}\vec q_{1}  - \sqrt{2(d-\Gamma)} w_2(l_2,t)\right]^2 \\
&\qquad -\frac{1}{2}  \Delta \vec q^{\rm T} Q  \vec q^{\rm T} -\int_0^{l_2} \kappa |z_x(-x,t)|^2 dx\\
&\qquad -\Gamma |w_2(l_2,t)|^2 - \xi_1 |w_1(l_2,t)|^2\le 0.
\end{array}
\end{eqnarray}
\end{lemma}
\begin{proof} By a direct calculation,
\begin{eqnarray*}
\begin{array}{ll}
&\frac{dE_{\rm{hybrid}}(t)}{dt}=\frac{1}{2}  ( \vec q^{\rm T}A^{\rm {T}} +  w_2(l_2,t) \vec b^{\rm T} )P\vec q \\
   & ~ + \frac{1}{2}   \vec q^{\rm T} P ( A \vec q+  w_2(l_2,t) \vec b)    +\int_0^{l_1} ( \kappa z {{z}}_{xx})(-x,t) dx\\
      &~ +\int_0^{l_2}\left[ w_1 (\alpha {u}_{1,x} -\gamma\beta  {u}_{2,x}) + w_2 (\beta {u}_{2,x} +\gamma \beta  {u}_{1,x})\right. \\
     &~ \left.+\alpha_1 u_1 {w}_{1,x}  +\beta (\gamma u_1-u_2) (\gamma {w}_{1,x}- {w}_{2,x})\right](x,t) dx, \\
            &=\frac{1}{2}   \vec q^{\rm T}(-\vec q_{1} \vec{q}_1^{\rm{T}}-\Delta Q)  \vec q  +   w_2(l_2,t)   \vec q^{\rm T} (\vec c+\sqrt{2(d-\Gamma)} \vec q_1)\\
            &~+2 w_2(l_2,t) \vec q^{\rm T} P \vec b- \xi_1 |w_1(l_2,t)|^2 \\
             &~ - \vec q^{T} (t) \vec c w_2(l_2,t)-d|w_2(l_2,t)|^2- \kappa\int_0^{l_1}{|{z}}_{x}(-x,t)|^2 dx.
       \end{array}
\end{eqnarray*}
Finally, by the assumptions in Section \ref{design}, \eqref{dissp2} follows.
\end{proof}
Define the Lyapunov function by
\begin{eqnarray}\label{Lyp2}
L_{\rm{hybrid}}(t)=E_{\rm{hybrid}}(t)+  \delta F(t)
\end{eqnarray}
where $F(t)$ is defined identically as in \eqref{Lyp1}.
\begin{lemma}\label{lem7}
Assume \eqref{constants}. Then, $F(t)$ in \eqref{Lyp2} satisfies the following estimate
\begin{eqnarray} \label{eqrfd3}
\begin{array}{ll}
\frac{dF(t)}{dt}&\le \frac{4b_1l_2}{l^2_1}\left(\vec q^{\rm{T}}(t)P\vec q(t)-2E(t)\right)\\
&\quad+\frac{a_1 l_2}{2}\left\{\left(\rho+2\frac{2\xi^2_1}{\alpha_1}\right)|w_1(l_2,t)|^2\right. \\
&\quad+\left(\mu+\frac{4\gamma^2 d^2}{\alpha_1}+\frac{2d^2} {\beta}\right)|w_2(l_2,t)|^2\\
&\quad \left.+\left(\frac{4\gamma^2}{\alpha_1}+\frac{2}{\beta}\right)\|\vec c\|^2\|\vec q(t)\|^2\right\}.
\end{array}
\end{eqnarray}
\end{lemma}
\begin{proof}
First, multiply both sides of the first  equation of \eqref{main2}  by $c_1 z(-x,t)$, and dot product the matrix equation of \eqref{main2}    by the vector of multipliers  \eqref{giant}. Integrate the results in $x$  from $0$ to $l_2$. After a couple of integration by parts and the boundary conditions \eqref{main-aug1}-\eqref{main-aug2}, we have the identity below, which is the counterpart of \eqref{mult1},
\begin{eqnarray}\label{mult11}
\begin{array}{ll}
& a_1\frac{d}{dt}\int^{l_2}_0  x\left[\rho  w_1u_1+\mu  u_2 w_2  \right] (x,t)dx\\
& =-\frac{a_1}{2}\int^{l_2}_0 \left[\rho |w_1|^2 +\mu  |u_2|^2 +\alpha_1 |u_1|^2 \right.\\
& \left. + \beta |\gamma u_1-u_2|^2 \right](x,t)dx+\frac{a_1 l_2}{2} \left[ \rho |w_1|^2 +\mu |w_2|^2 \right. \\
& \left.+\frac{\left|\gamma \left(\vec c^{\rm {T}} \vec q+dw_2\right)+ \xi_1 w_1\right|^2}{\alpha_1}~ +\frac{\left|\vec c^{\rm {T}} \vec q+dw_2\right|^2}{\beta}\right](l_2,t).
\end{array}
\end{eqnarray}
By the  Cauchy-Schwarz's  and  Wirtinger's inequalities,
\begin{eqnarray*}
\begin{array}{ll}
&\frac{dF}{dt}\le \frac{-4b_1 l_2 \kappa}{l^2_1}\int^{l_1}_0 |z(-x,t)|^2 dx\\
&\quad + 2 \int^{l_2}_0 \frac{b_1}{2}\left[\left(\sqrt{\frac{\alpha_1}{\rho}}+2\sqrt{\frac{\beta}{\mu}}\frac{\mu \gamma^2}{\rho}\right)\rho |w_1|^2+\sqrt{\frac{\alpha_1}{\rho}}\alpha_1 |u_1|^2\right.\\
&\quad \left. +2\sqrt{\frac{\beta}{\mu}}\mu |w_2|^2 +\sqrt{\frac{\beta}{\mu}}\beta |\gamma u_1 -u_2|^2 \right](x,t)dx\\
&\quad -\frac{a_1}{2}\int^{l_2}_0 \left[\rho |w_1|^2+\mu |w_2|^2 \right. \\
&\left.\quad +\alpha_1 |u_1|^2+\beta \left|\gamma u_1-u_2\right|^2\right](x,t)dx\\
&\quad +\frac{a_1 l_2}{2}\left(\rho |w_1(l_2,t)|^2+\mu |w_2(l_2,t)|^2\right.\\
&\quad +\frac{2\gamma^2\left|\vec c^{\rm {T}} \vec q+dw_2(l_2,t) \right|^2 + 2\xi_1^2 |w_1(l_2,t)|^2}{\alpha_1}\\
&\left.\quad +\frac{2 \|\vec c\|^2 \|\vec q\|^2+ d^2 |w_2(l,t)|^2}{\beta}\right).
\end{array}
\end{eqnarray*}
Finally, by  the  Cauchy-Schwarz's inequality and \eqref{constants}, \eqref{eqrfd3} follows.
\end{proof}

  	\begin{theorem} Assume \eqref{constants}. Suppose that there exists a constant
	\begin{eqnarray}\label{del-as}
	\begin{array}{ll}
	0<\delta <\frac{1}{a_1l_2}{\rm {min}} \left(\frac{a_1l_2}{M},\frac{2\xi_1}{\rho+\xi_1^2\frac{2}{\alpha_1}},\right.\\
	\left.\frac{2\Gamma}{\mu+2d^2\frac{\alpha+\gamma^2\beta}{\alpha_1\beta}},\frac{\Delta \vec q^{\rm{T}} Q\vec q}{\frac{2\alpha+2\gamma^2\beta}{\alpha_1\beta}\|\vec c\|^2\|\vec q(t)\|^2+\frac{8a_1b_1 l_2 k_2\kappa}{l_1^2}\vec q^{\rm{T}}P\vec q} \right).
	\end{array}
	\end{eqnarray}
Then, for all initial conditions, the energy $E_{\rm {hybrid}}(t)$ along the solutions of  \eqref{main2}  with \eqref{main-aug1}-\eqref{main-aug2}  decays exponentially, i.e.
\begin{eqnarray*}
E_{\rm {hybrid}}(t)\le\frac{1+M\delta}{1-M\delta} e^{-\sigma t} E_{\rm {hybrid}}(0), \quad \forall t>0
\end{eqnarray*}
where $\sigma:=\delta (1-M\delta) \frac{8 b_1 l_2\kappa}{l_1^2}.$
\end{theorem}
\begin{proof} By \eqref{MKY},  \eqref{del-as}, and Lemmas \ref{lem6} and \ref{lem7},
\begin{eqnarray*}
\begin{array}{ll}
 \frac{dL_{\rm{hybrid}}(t)}{dt}\le -\kappa \int^{l_1}_0 |z_x(-x,t)|^2 dx -\frac{8b_1l_2\delta }{l^2_1}E(t)\\
   \quad+\frac{4b_1l_2\delta }{l^2_1}\vec q^{\rm{T}}(t)P\vec q(t)  +\frac{a_1l_2}{2} \delta\left(\frac{2\alpha_1+4\gamma^2\beta}{\alpha_1\beta}\right)\|\vec c\|^2\|\vec q(t)\|^2  \\
   \quad -\frac{1}{2}  \Delta \vec q^{\rm T} Q  \vec q^{\rm T}-\left(\xi_1-\frac{\delta a_1 l_2}{2}\left(\rho+2\frac{2\xi^2_1}{\alpha_1}\right)|w_1(l_2,t)|^2\right)\\
   \quad -\left(\Gamma-\frac{\delta a_1 l_2}{2}\left(\mu+\frac{4\gamma^2 d^2}{\alpha_1}+\frac{2d^2} {\beta}\right)|w_2(l_2,t)|^2\right)\\
    \quad-\frac{1}{2}  \left[ \vec q^{\rm T}\vec q_{1}  - \sqrt{2(d-\Gamma)} w_2(l_2,t)\right]^2 \le  -\frac{8b_1l_2\delta }{l^2_1}E(t).
\end{array}
\end{eqnarray*}
Note that since $E(t)\le E_{\rm{hybrid}}(t),$ Lemma \ref{lem2} is still valid, i.e.  for $0<\delta < \frac{1}{M},$  $L(t)$ is equivalent to $E_{\rm{hybrid}}(t),$ i.e.,
$$(1-M \delta) E_{\rm{hybrid}}(t)\le L(t)\le (1+M\delta)E_{\rm{hybrid}}(t).$$
Therefore, by the Gr\"onwall's inequality and \cref{lem7},
\begin{eqnarray*}
\begin{array}{ll}
L_{\rm{hybrid}}(t)\le L_{\rm{hybrid}}(0)e^{-\sigma t} \le (1+M\delta) E_{\rm{hybrid}}(0) e^{-\sigma t}.
\end{array}
\end{eqnarray*}
Hence, \eqref{del-as} follows.
\end{proof}
\begin{remark}
Note that $\sigma(\delta(\xi_1,\vec b, \vec c, d))$ makes its maximal value
$\sigma_{max}(\delta)= \frac{2b_1 l_2 \kappa  }{l_1^2 M}$ achieved at  $\delta=\frac{1}{2M},$
similar to the result in Theorem \ref{mainthm1}. However, finding optimal values $\xi_1,\vec b, \vec c, d, {\bf A}$  for achieving  $\delta=\frac{1}{2M}$ is  a more complicated task.
\end{remark}

\section{Model Reduction by Finite Differences}
\label{Simu}

For simplicity, assume that $l_1=l_2=L.$ Let $N\in\N$ be given. Defining the mesh size $h:=\frac{L}{N+1}$, the following discretization of the interval $[0,L]$ is considered
$0=x_0<x_1<...<x_{N-1}<x_N<x_{N+1}=L.$
For this method, we also consider the middle points of each subinterval denoted by $\{x_{j+\frac{1}{2}}\}_{j=0}^{N+1},$ where i.e. $x_{j+\frac{1}{2}}=\frac{x_{j+1}+x_j}{2}.$ 	Define following finite average and difference operators $
		u_{j+\frac{1}{2}}:=\frac{u_{j+1}+u_j}{2}, \quad \delta_xu_{j+\frac{1}{2}}:=\frac{u_{j+1}-u_j}{h}.$
Note that considering the odd-number derivative at the in-between nodes  provides  higher-order approximations.

For each $j=1,2,\ldots,N,$ define $z_{-j}(t):=z(-x_j,t),$
$(u^1, u^2, w^1, w^2)_j:=(u_1, u_2,  w_1,  w_2)(x_j), $
Moreover, they can be represented in two separate ways by the preceding  and following nodes \cite{Ren}.   Letting
${{\tilde A}}_{4\times 4}:=\begin{bmatrix}
		0& I\bigotimes \delta_x \\
			A\bigotimes \delta_x& 0      \\
			\end{bmatrix}, $  the order-reduced Finite-Difference formulation of the hybrid model \eqref{main2} with \eqref{main-aug1}-\eqref{main-aug2} is proposed as the following
\begin{eqnarray}
\begin{array}{ll}
\label{main3}
		\left\{
	 \begin{array}{ll}
	z_{-j,t}(t)-\delta_x^2 z_{-j}(t)=0,\\
\tilde M\begin{bmatrix}
		 u^1_{j{+\frac{1}{2}},t}\\
	          u^2_{j{+\frac{1}{2}},t}\\
		w^1_{j{+\frac{1}{2}},t} \\
		w^2_{j{+\frac{1}{2}},t}     \\
			\end{bmatrix}-
{ {\tilde A}} \begin{bmatrix}
		 u^1_{j{+\frac{1}{2}}}\\
	          u^2_{j{+\frac{1}{2}}}\\
		w^1_{j{+\frac{1}{2}}} \\
		w^2_{j{+\frac{1}{2}}}     \\
			\end{bmatrix}=0,\\
			\vec q_t(t)= {\bf A} \vec q(t) + \vec b w^2_{N+1} (t), \qquad  t\in\mathbb{R}^+,
			\end{array}
			\right.
			\end{array}
			\end{eqnarray}
			\begin{eqnarray}
			\left\{ \begin{array}{ll}
		z_{-(N+1)}(t)=0, \quad z_0(t) = w^1_0 (t), ~~w^2_0(t)=0, \\
		\delta_x z_0(t)=\frac{z_0(t)-z_{1}(t)}{h}=\alpha u^1_0 (t)-\gamma\beta u^2_0(t),&\\
A
 \begin{bmatrix}
				u^1_{N+1} \\
				u^2_{N+1}   \\
			\end{bmatrix}(t) = \begin{bmatrix}
				-\xi_1 w^1_{N+1}\\
				- \vec c^{\rm{T}} \vec q -d w^2_{N+1}  \\
			\end{bmatrix}(t),~  t\in\mathbb{R}^+,\\
			z_{-j}(0)=z^0(-x_j), \quad \vec q(0)=\vec\eta,\\
\left(u^1,u^2,w^1, w^2 \right)_j(0)\\
\qquad =(v^0_x, p^0_x, v^1,p^1)(x_j),~j=0,1,\ldots, N+1.
		\end{array}\right.
	\end{eqnarray}
with the choice of the discretized energy $E_{\rm {hybrid}}(h,t)$
\begin{eqnarray*}
\begin{array}{ll}
	E_{\rm {hybrid}}(h,t):=\frac{1}{2} \vec q^{\rm T}(t) P \vec q(t)\\
	+\frac{h}{2} \sum\limits_{j=0}^N\left\{ |z_{-j}(t)|^2+\left(A   \begin{bmatrix}
				u^1_{j+\frac{1}{2}}(t)  \\
				u^2_{j+\frac{1}{2}}(t)      \\
			\end{bmatrix}\right) \cdot \begin{bmatrix}
				u^1_{j+\frac{1}{2}}(t) \\
				u^2_{j+\frac{1}{2}}(t)       \\
			\end{bmatrix}  \right.\\
\left.			
	+\left( M\begin{bmatrix}
				w^1_{j+\frac{1}{2}}(t)  \\
				w^2_{j+\frac{1}{2}} (t)   \\
			\end{bmatrix}\right) \cdot \begin{bmatrix}
				w^1_{j+\frac{1}{2}}(t)\\
				w^2_{j+\frac{1}{2}} (t)  \\
			\end{bmatrix}		 \right\} dx.
\end{array}\end{eqnarray*}
which is dissipative along \eqref{main3}. Note that the discretized model above can be  adjusted for the static feedback case.

\section{Conclusion\& Future Work}
The immediate future work is to analyze the  model reduction of the original coupled PDEs \eqref{main}, proposed in Section \ref{Simu}. The power of this type of semi-discretization lies on the fact that the Lyapunov approach and the multipliers used in Sections \ref{sec2} and \ref{sec3} can be adapted in a discrete fashion to prove the uniformly exponential stability of \eqref{main3}, uniformly as the discretization parameter $h$ tends to zero \cite{Jin,Ren,Zheng}.

The system \eqref{main} also models the fluid-structure interactions for the linearized one-dimensional Navier-Stokes equations (heat equation) and the solid structure (piezoelectric beam). Therefore,  the proposed methodology can also solve the same feedback stabilization problem \cite{Zua1}. A potential extension of this work is considered to be transmission line settings with thermoelastic, elastic, and piezoelectric beams, as in \cite{Serge}.

\end{document}